\documentclass{commat}

\usepackage{float}

\title{%
    Relaxation in one-dimensional tropical sandpile
    }

\author{%
   Mikhail Shkolnikov
    }

\authorinfo{Institute of Mathematics and Informatics, Bulgarian Academy of Sciences}{mikhail.shkolnikov@hotmail.com}

\abstract{%
A relaxation in the tropical sandpile model is a process of deforming a tropical hypersurface towards a finite collection of points. We show that, in the one-dimensional case, a relaxation terminates after a finite number of steps. We present experimental evidence suggesting that the number of such steps obeys a power law. 
    }



\VOLUME{31}
\NUMBER{3}
\firstpage{21}
\DOI{https://doi.org/10.46298/cm.10483}

\begin{document}

\section{Introduction}

The sandpile model was discovered independently several times and in different contexts (see \cite{LP10}). It became especially popular when it was proposed as a prototype for self-organized criticality \cite{BTW87}. This somewhat vague concept can be defined in various complementing ways, the most straightforward is that the system has no tuning parameters and demonstrates power-laws. We describe a very simple model (see Figure \ref{relax09}) having such property. 

Until very recently \cite{K21}, the tropical sandpile model has been discussed only in two-dimensional case. It arises as a scaling limit of the original sandpile model in the vicinity of the maximal stable state (this is formally stated in the case of lattice polygonal domains in \cite{KS16} and proven for general convex domains in \cite{KS15}) and was studied numerically in \cite{KL18}, where it was shown to exhibit a power law providing the first example of a continuous self-organized criticality. The later direction is further explored in \cite{KP23}.

The setup for the tropical sandpile model is as follows. Consider a compact convex domain $\Omega\subset\mathbb{R}^d.$ A function $F\colon\Omega\rightarrow [0,\infty)$ is called an $\Omega$-tropical series if it vanishes on $\partial\Omega$ and can be presented as \[F(z)=\inf_{v\in\mathbb{Z}^d} (a_v+z\cdot v).\] The numbers $a_v\in\mathbb{R}$ are called the coefficients of $F.$ The coefficients of $F$ are not uniquely defined. However, there is a canonical choice, i.e. we set them to be as minimal as possible. 

\begin{figure}
    \centering
    \includegraphics[width=.5\textwidth]{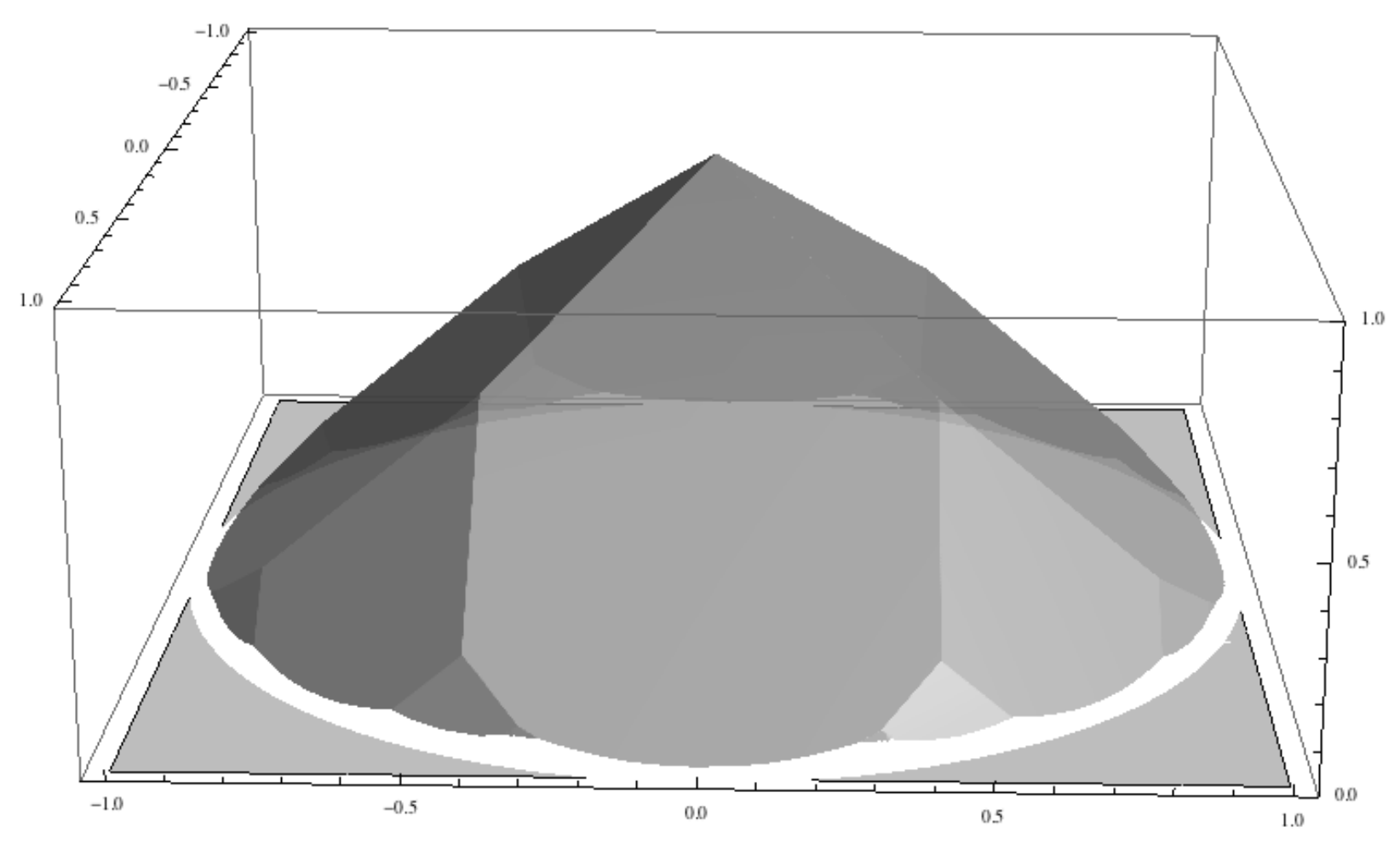}
    \caption{A plot of the tropical series $z\mapsto\inf_{v\in\mathbb{Z}^2\backslash\{0\}} (|v|+z\cdot v)$ on the unit disk.}
    \label{2dSeries}
\end{figure}

For example, take $\Omega$ to be a disk $\{z\in\mathbb{R}^2:|z|\leq 1\}.$ Then, $\inf_{v\in\mathbb{Z}^2\backslash\{0\}} (|v|+z\cdot v)$ is an $\Omega$-tropical series (see Figure \ref{2dSeries} for the plot of this function). We see that the ``monomial'' corresponding to $0\in\mathbb{Z}^2$ doesn't participate in the formula, but in the canonical choice of the coefficients we need to take $a_0$ to be $1,$ that is the maximal value of the series, which is attained at the origin.

The initial state of the model $0_\Omega$ is an $\Omega$-tropical series vanishing on the whole $\Omega.$ Its coefficient corresponding to $0\in\mathbb{Z}^d$ is $0$ and, in the canonical form, its coefficient for $v\in\mathbb{Z}^d\backslash\{0\}$ is $-\min_{z\in\Omega} z\cdot v.$ 

For a point $p\in\Omega^\circ,$ we define an idempotent operator $G_p$ acting on the space of $\Omega$-tropical series. In short, $G_pF$ is the result of increasing at most one coefficient in the canonical form of $F,$ so that it attains a break (i.e. becomes not smooth) at $p.$  More explicitly, if $F$ is linear at $p,$ then there exist a unique $w\in\mathbb{Z}^d$ such that $F(z)=a_w+z\cdot w$ for $z$ in a neighborhood of $p,$ and we take $G_pF$ to be $z\in\Omega\mapsto \inf_{v\in\mathbb{Z}^d} (b_v+z\cdot v),$ where $b_v=a_v$ for $v\in\mathbb{Z}^d\backslash\{w\},$\ \[b_w=\min_{v\in\mathbb{Z}^d\backslash\{w\}}(a_v+p\cdot v)-p\cdot w\] and $a_v$ are the canonical coefficients of $F.$ See Figure \ref{Gp1d} for a one-dimensional example and Figure \ref{shrinkPhi} for a two-dimensional geometric version. Observe, also, that the series plotted on Figure \ref{2dSeries} is $G_p0_\Omega$ for $\Omega$ equal to the unit disk and $p$ equal to its center.

The operator $G_p$ is a tropical counterpart of adding a grain at $p$, relaxing and then removing the grain in the sandpile model, when one works in the tropical sector, where all states are made of sandpile solitons. In its infinitesimal form, i.e. that increasing a coefficient of a tropical polynomial corresponds to sending a wave, this statement appears as Corollary 2.13 in \cite{KS20}, in the planar case, and as Proposition 3.5 \cite{K23}, for arbitrary dimension. For its finite form, i.e. after applying enough waves, see, for example ``{\it A sketch of a proof}'' of Lemma 1 in \cite{KL18}.

Consider a collection of points $p_1,\dots,p_n\in\Omega^\circ.$ A relaxation is a sequence of the form: \begin{equation} F_m=G_{p_{k_m}}G_{p_{k_{m-1}}}\dots G_{p_{k_1}}0_\Omega\tag{1}\end{equation} of $\Omega$-tropical series, where $k_1,k_2,\dots\in\{1,\dots,n\}$ is a sequence of indices taking each value infinitely many times. For $d=2,$ it was shown in \cite{KS18} that $F_m$ uniformly converges to $G_{\{p_1,\dots,p_n\}}0_\Omega,$ the minimal $\Omega$-tropical series not smooth at $p_1,\dots,p_n.$ In fact, the argument works equally well for all $d$ (see \cite{K21}).

However, unless $\Omega$ is a lattice polytope and $p_1,\dots,p_n\in\mathbb{Z}^d,$ it is not clear if a relaxation terminates after a finite number of steps. We prove the following.

\begin{theorem}
If $d=1,$ then the sequence $F_m$ stabilizes.
\end{theorem} 
It is reasonable now to consider a question:   {\it What is the distribution for the length of relaxation?}

To make this question more precise, for points $p_1,\dots,p_n$ we define the length of relaxation $L(p_1,\dots, p_n)$ as the minimal number $N$ such that \[(G_{p_n}\dots G_{p_1})^N 0_\Omega=G_{\{p_1,\dots,p_n\}}0_\Omega.\] We want to look at the distribution of $L(p_1,\dots, p_n)$ when $p_1,\dots ,p_n\in\Omega^\circ$ are taken as independent uniform random variables. Our computer simulation suggests the presence of power-laws (see Figure \ref{plots}), surprisingly, already for $n=2$. 

\section{Stabilization}

\begin{figure}
    \centering
    \includegraphics[width=.66\textwidth]{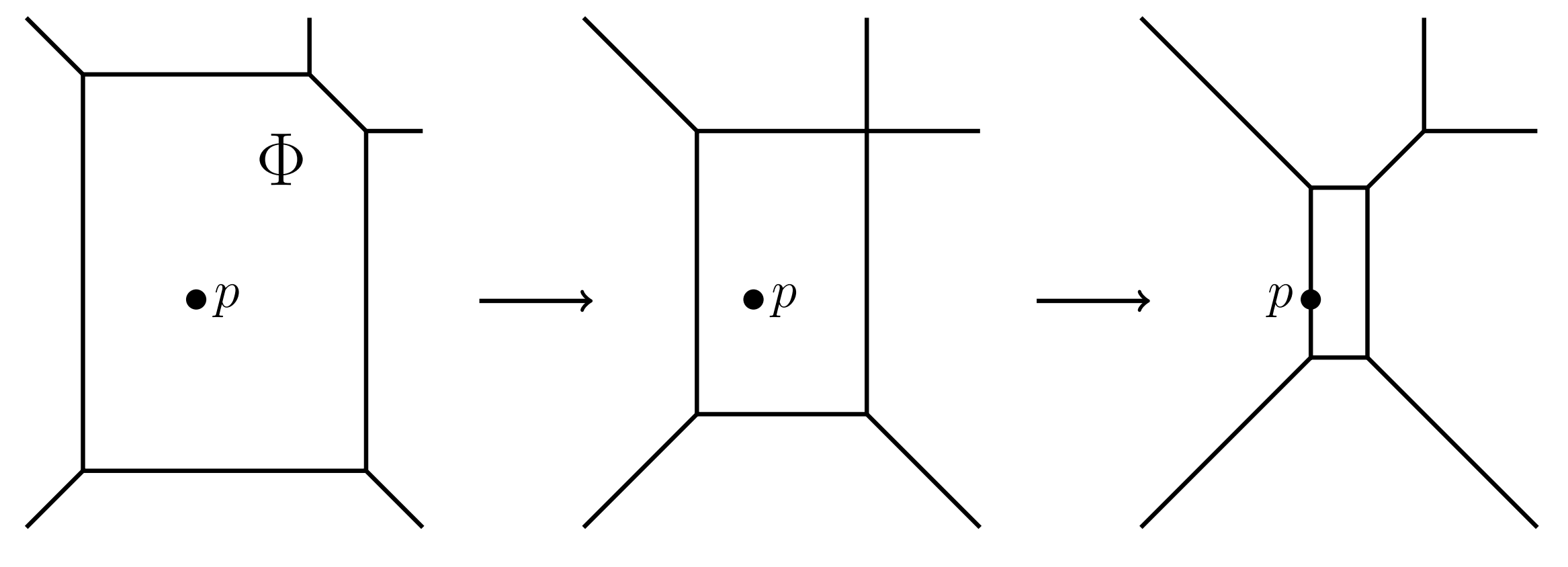}
    \caption{Two-dimensional illustration of the geometric interpretation of the operator $G_p$ seen as a result of a continuous process that shrinks a face to which the point $p$ belongs. }
    \label{shrinkPhi}
\end{figure}

The operator $G_p$ has a nice geometric interpretation in terms of hypersurfaces. An $\Omega$-tropical series $F$ defines its $\Omega$-tropical hypersurface $H$ as a locus of all points $z\in\Omega^\circ$ where $F$ is not smooth. If $p\in H$ then $G_p F=F.$ Otherwise, the hypersurface defined by $G_p F$ may be thought as the result of shrinking the connected component of $\Omega^\circ\backslash H$ containing $p$ (see Figure \ref{shrinkPhi}).  We will describe explicitly how this works in the one-dimensional case.

Let $\Omega$ be an interval. A hypersurface defined by an $\Omega$-tropical series $F$ is just a discrete set of points $H\subset\Omega^\circ$ over which the graph of $F$ breaks. We incorporate multiplicities $\mu\colon H\rightarrow \mathbb{Z}_{\geq 1}$ for these points by computing the second derivative, i.e. \[\frac{d^2}{dx^2}F(x)=-\sum_{h\in H}\mu(h)\delta(x-h),\tag{2}\] where  $\delta$ is the Dirac delta function. In other words, if $p\in H$ then $\mu(p)$ is equal to the difference between the slopes of the linear pieces of $F$ to the left and to the right from $p$.

In the rest of this note, we assume that $H$ is {\it finite}, i.e., $F$ is the restriction to $\Omega$ of a tropical polynomial vanishing on $\partial\Omega$. We call such $F$ an $\Omega$-tropical polynomial.

One can restore $F$ from $H$ and $\mu$ (note that there is no constant and linear term ambiguity since $F$ has to vanish at the boundary of the interval $\Omega$). However, not every finite collection of points with multiplicities is defined by an $\Omega$-tropical polynomial. Indeed, performing twice an indefinite integration of the right-hand side of (2), we get a two-dimensional space of functions of the form $F_{\alpha,\beta}(x)=f(x)+\alpha x+\beta,$ where $f$ is a piece-wise linear function with integral slopes and $\alpha,\beta$ are any real numbers.  There is a unique choice of $\alpha$ and $\beta$ such that $F=F_{\alpha,\beta}$ vanishes on $\partial\Omega.$ Unless $\alpha$ is an integer, $F$ fails to be an $\Omega$-tropical polynomial. We will use the following criterion.

\begin{figure}
    \centering
    \includegraphics[width=0.5\textwidth]{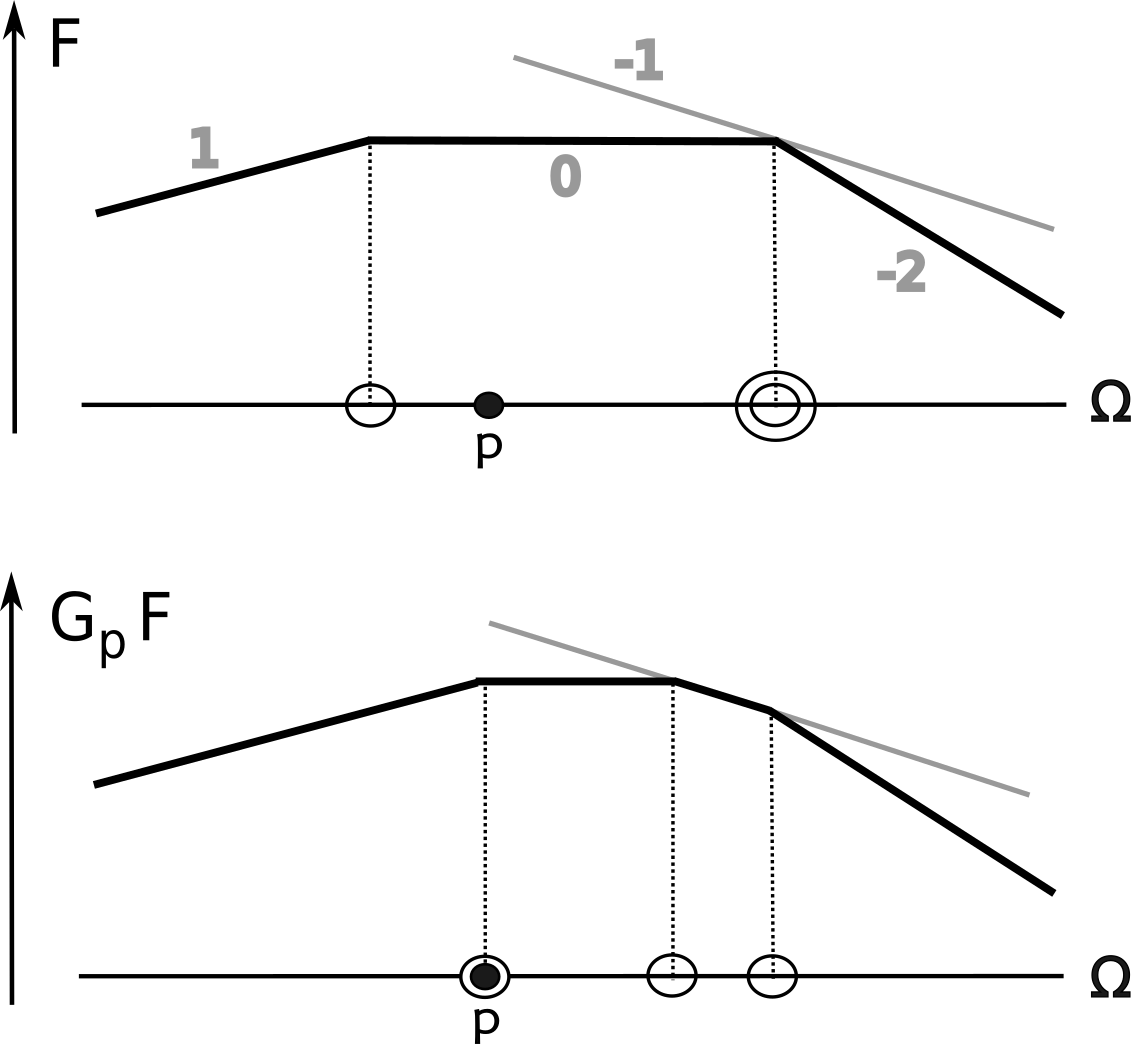}
    \caption{One-dimensional illustration of the operator $G_p.$ Plots of $F$ and $G_pF$ are shown in bold, the slopes of linear pieces are written near them in gray. In this example, the original hypersurface $H$ defined by $F$  contains a point of multiplicity two, since there is an extra monomial (its plot is shown in gray) in the canonical form of $F$ contributing at this point, or, equivalently, since the slope jumps by two at this point. Applying $G_p$ at the level of tropical polynomials means increasing the coefficient of the constant term (moving upwards the horizontal segment of the plot in this example) until the function has a break at $p.$ At the level of sets of non-linearity of tropical polynomials, $G_p$ moves two points towards $p$ until at least one of them coincides with $p.$}
    \label{Gp1d}
\end{figure}

\begin{proposition} Let $\Omega=[0,1]$. A finite set  $H\subset (0,1)$ with multiplicities $\mu$ is defined by an $\Omega$-tropical polynomial if and only if $\sum_{h\in H}\mu(h)h$ is an integer. 
\end{proposition}
\begin{proof} For $h\in(0,1)$ let $f_h(x)$ be a definite double integral of $-\delta(x-h),$ i.e., \[f_h(x)=-\int_0^x\int_0^t\delta(s-h)dsdt.\] Note that $f_h(x)=\min(0,h-x).$ Therefore, its value at $0$ is $0$ and at $1$ is $h-1.$ To make $\alpha x+\sum_{h\in H}\mu(h)f_h(x)$ vanish at $1$ we should take $\alpha=-(\sum_{h\in H}(h-1))x,$ which is an integer if and only if $\sum_{h\in H}h$ is an integer.
\end{proof}

To express $G_p$ in a closed form, it will be convenient to encode $F$ by a function \[M_F\colon\Omega\rightarrow\mathbb{Z}_{\geq 0}\cup\{\infty\}\] defined as $M_F(\partial\Omega)=\{\infty\},$ $M_F(\Omega^\circ\backslash H)=\{0\}$ and $M_F|_H=\mu.$ Assume $p$ belongs to a connected component $(a,b)$ of  the complement of $H$ in $\Omega^\circ$. Then \[M_{G_pF}(x)=M_F(x)-\delta_{a,x}-\delta_{b,x}+\delta_{a+c,x}+\delta_{b-c,x},\] where $\delta_{\cdot,\cdot}$ is the Kronecker delta and $c=\min(p-a,b-p)$. In plain words, $G_p$ moves by $c$ the ends of the connected component towards $p,$ see Figure \ref{Gp1d}. 

\begin{remark} $G_p$ doesn't produce points with multiplicities greater than $2.$
\end{remark}

For example, let $\Omega=[a,b]$ and $p=p_1\in (a,b).$ If $2p\neq a+b,$ then the set of points defined by $G_p 0_\Omega$ is $\{p, a+b-p\},$ and the multiplicity of each point is $1.$ If $2p=a+b,$ then the set consists of a single point $p$ with multiplicity $2.$ We see that for one point, the relaxation terminates after one step.

\begin{figure}
    \centering
    \includegraphics[width=0.5\textwidth]{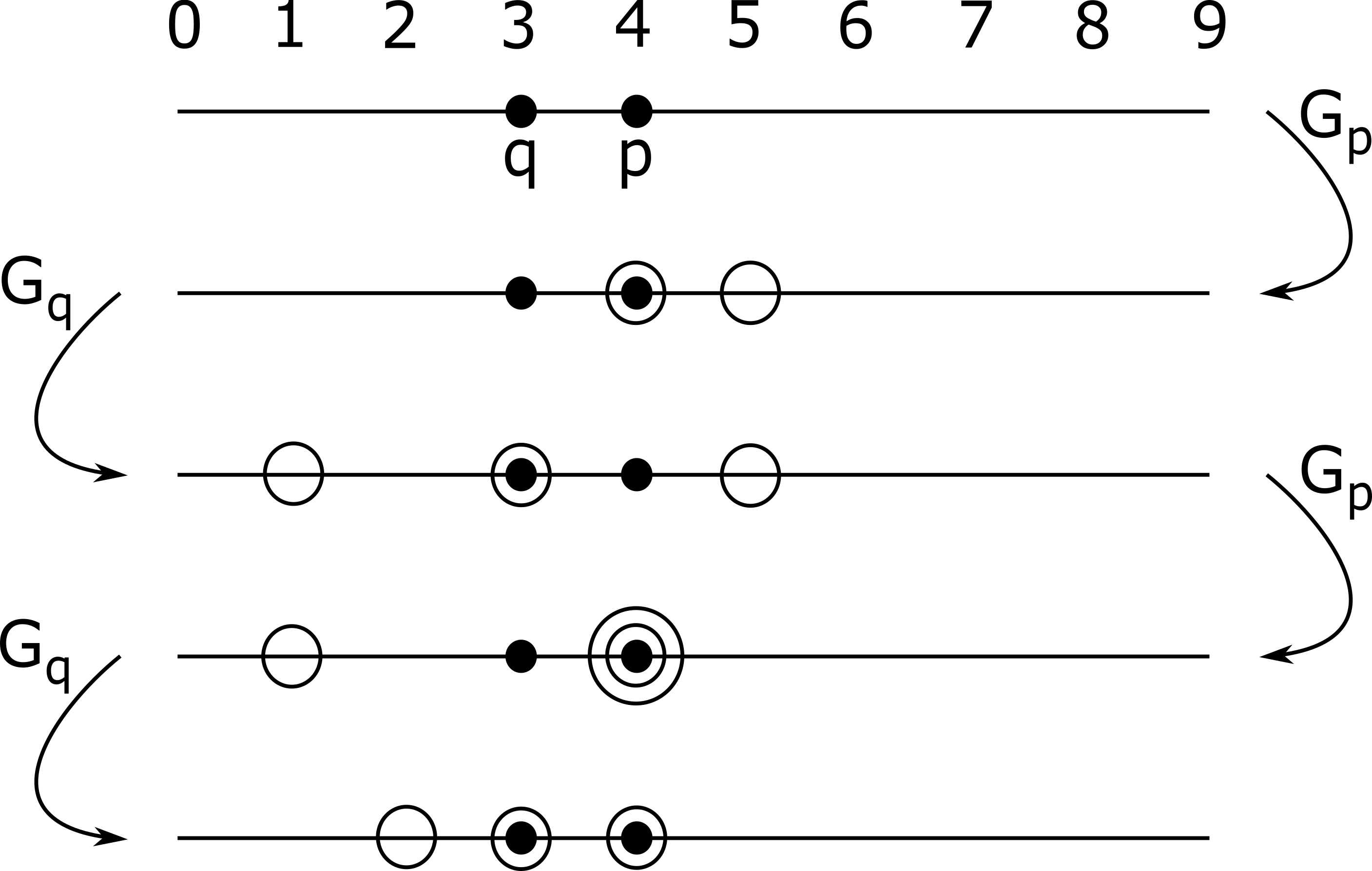}
    \caption{Relaxation for points $p=4$ and $q=3$ on $\Omega=[0,9],$ represented by hypersurfaces (denoted by collections of circles -- two concentric circles circle mean a point of multiplicity two) defined by $\Omega$-tropical polynomials $0_\Omega,$ $G_p 0_\Omega,$ $G_qG_p 0_\Omega,$ $G_p G_qG_p 0_\Omega,$ etc. The lowest picture is stable, i.e. the points $p$ and $q$ belong to the hypersurface defined by the result of the relaxation.}
\label{relax09}
\end{figure}

For a less trivial and more concrete example of a relaxation process, take $\Omega=[0,9],$ $p=p_1=4,$ and $q=p_2=3$.  Then, $G_p 0_\Omega$ defines points $4$ and $5;$ $G_qG_p 0_\Omega$ defines points $1,\ 3$ and $5;$  $G_pG_qG_p 0_\Omega$ defines $1$ with multiplicity $1$ and $4$ with multiplicity $2;$ finally, $G_qG_pG_qG_p 0_\Omega$ defines $2,\ 3,$ and $4$ (see Figure \ref{relax09}).

We now proceed to proof of the stabilization Theorem.

\begin{proof}
We describe first $F_\infty=G_{\{p_1,\dots,p_n\}}0_\Omega,$ the limit of $F_m$ defined by (1). Without loss of generality, assume that $\Omega$ is $[0,1]$ and that the points $p_1,p_2\dots p_n\in (0,1)$ are distinct.  Let $q\in[0,1)$ be the fractional part of $-\sum_{j=1}^np_j.$

\begin{lemma} The structure of the set $H_\infty$ with multiplicities $\mu_\infty$ defined by $F_\infty$ depends on the position of  $q$:
\begin{itemize}
\item if $q=0,$ then $H_\infty=\{p_1,\dots,p_n\}$ and all multiplicities are $1$;
\item if there exist $j\in\{1,\dots,n\}$ such that $q=p_j,$ then $H_\infty=\{p_1,\dots,p_n\}$ and all multiplicities are $1$ except for $\mu_\infty(p_j)=2$;
\item otherwise, $H_\infty=\{q, p_1,\dots,p_n\}$ and all multiplicities are $1.$
\end{itemize} 
\end{lemma}

\begin{proof}
Notice that the number of points in $(H_\infty,\mu_\infty)$ counted with multiplicities is equal to the difference of slopes of $F_\infty$ at $0$ and $1,$ and the absolute values of these slopes are minimized by $F_\infty$ in the class of $\Omega$-tropical polynomials not smooth at $p_1,\dots,p_m,$ since $F_\infty$ is the pointwise minimum of all such functions (Definition 5.3 and \cite[Proposition~6.1]{KS18}) Therefore, $(H_\infty,\mu_\infty)$ is determined by the condition that it is the smallest multi-set containing all $p_i$ and satisfying the criterion of the Proposition.
\end{proof}

Consider the third (generic) case when $q\neq 0$ and $q\neq p_j$ for all $j$. The convergence of $F_m$ to $F_\infty$ implies that the set $H_m$ defined by $F_m$ converges to $H_\infty$. Take $\varepsilon>0$ to be smaller than the half of a minimal distance between two points of $H_\infty.$ There exists $m_\varepsilon$ such that for all $m\geq m_\varepsilon$  the $\varepsilon$-neighborhood of every point in $H_\infty$ contains a unique point of $H_m,$ and vice versa. Let $p_{m,i}\in H_m$ be the point in the $\varepsilon$-neighborhood of $p_i.$ 

Denote by $P_l$ the set of all $p_i$ smaller than $q$ and by $P_r$ the set of all $p_i$ greater than $q.$ We prove the stabilization of relaxation separately for $P_l$ and $P_r,$ the proofs are identical. 

Let $p_s$ be the smallest element of $P_l.$ Note that $p_{m_\varepsilon,s}$ cannot be greater than $p_s$ since otherwise at some further step $m_s>m_\varepsilon$ of the relaxation, when applying $G_{p_s},$ we would increase the number of points in $H_{m_s}$ as compared with $H_{m_s-1}$. Therefore,  $p_{m,s}=p_s$ for $m\geq m_s.$ This implies that for the second smallest point $p_k$ in $P_l$ we have $p_k\geq p_{m_k,k},$ otherwise, applying $G_{p_k}$ at some further step $m_k>m_s$ would violate $p_{m_k,s}=p_s.$ Thus, $p_{m,k}=p_k$ for $m\geq m_k.$ Etcetera.

Going from smaller $p_i$ to greater ones we have a chain of stabilizations at points of $P_l$. This chain is interrupted by the point of $H_m$ in the neighborhood of $q,$ so we need to launch another chain of stabilizations over $P_r,$ going from greater to smaller points. 

In the first case of the Lemma, we don't have this effect, so we need to do a single chain. In the second case, we proceed as in the third case and prove the stabilization at $p_j=q$ after we work out all other points (just before the last step $m_{\text{last}}$ the point $p_j$ is between two nearby points of $H_{m_{\text{last}}-1}$). 
\end{proof}

A similar argument should work in all dimensions. Instead of one or two linear chains of stabilizations, for a generic configuration of points $p_1,\dots,p_n$, there might be several tree-like chains.  However, a special care is needed for non-generic configurations when cycles in these chains may appear.

\section{Length of relaxation}

In this section, we will touch on the behavior of $L(p_1,\dots,p_n)$ defined at the end of the introduction. Specific choice of a segment $\Omega$ is irrelevant (one can apply an affine reparameterization); therefore, we restrict our attention to $\Omega=[0,1].$ 

First we note that there is an obvious symmetry  \[L(p_1,\dots,p_n)=L(1-p_1,\dots,1-p_n).\tag{3}\] On the other hand,  $L$ is sensitive to the order of its arguments. The closures of loci $L(p_1,\dots,p_n)=\text{const}$ are non-empty polytopal complexes with rational slopes.

For $n=1,$ the situation is very simple, i.e. $L(p)=1$ for all $p\in(0,1).$ For $n=2,$ we derive the following pictures. 

\begin{figure}[H]
    \centering
    \includegraphics[width=0.32\textwidth]{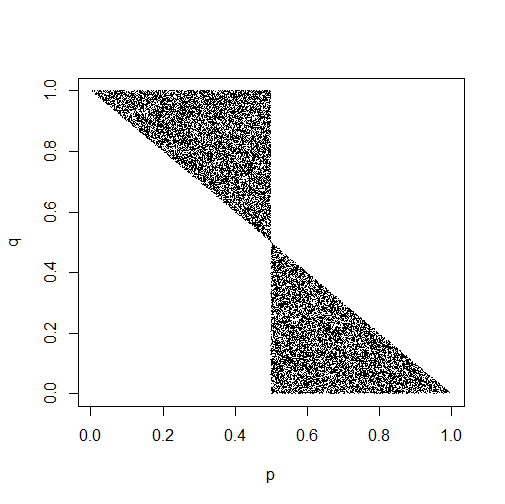}
    \includegraphics[width=0.32\textwidth]{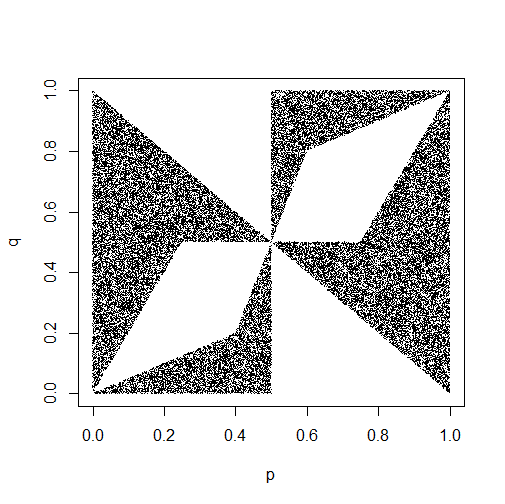}
    \includegraphics[width=0.32\textwidth]{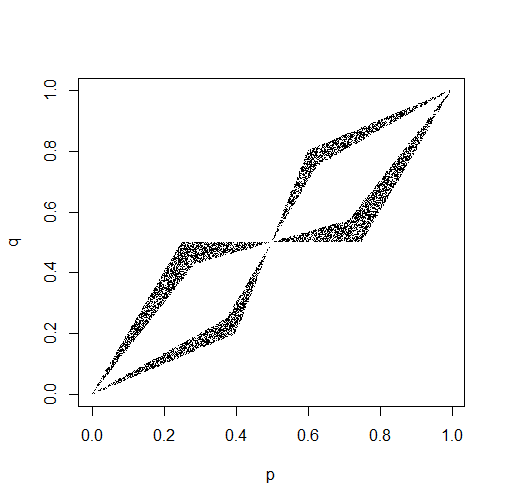}
    \includegraphics[width=0.32\textwidth]{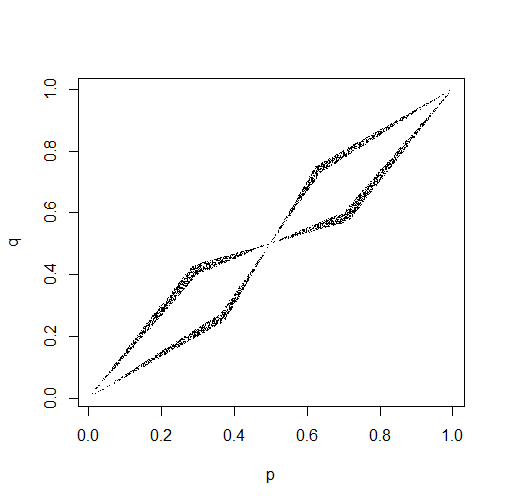}
    \includegraphics[width=0.32\textwidth]{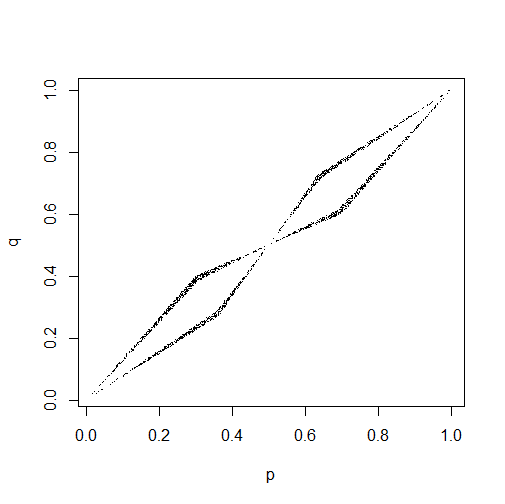}
    \includegraphics[width=0.32\textwidth]{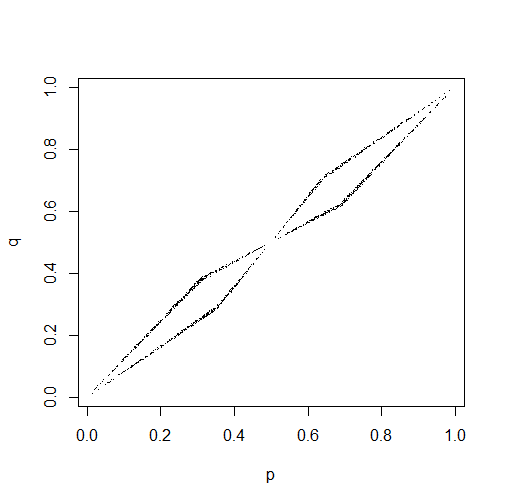}
    \caption{Numerical approximations for loci of $(p,q)\in(0,1)^2$ with length of relaxation $L(p,q)$ equal to $1,2,3,4,5$ or $6$ (left to right, up to down). The pictures are made in R.}
    \label{lenrel123456}
\end{figure}

\begin{figure}
    \centering
    \includegraphics[width=.67\textwidth]{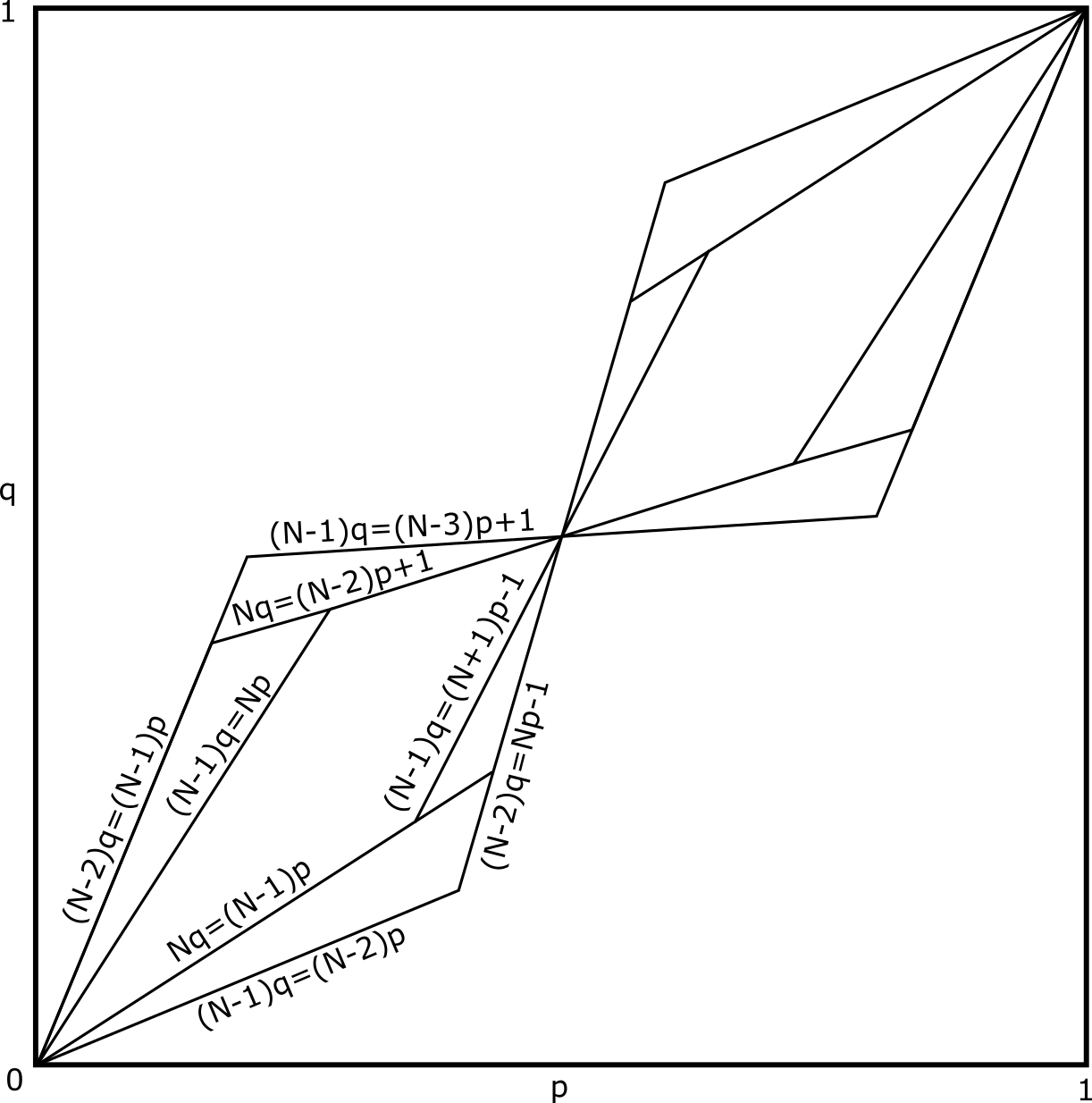}
    \label{lenrelN}
    \caption{For $N\geq 2,$ the closure of the locus $L(p,q)=N$ consists of eight triangles bounded by lines with equations written near them. The author has found this picture by carefully considering various cases for the positions of the points $p$ and $q$; an interested reader is encouraged to reproduce this picture by themselves. Hint: imagine first, how a relaxation goes if two points $p$ and $q$ are close to each other -- each intermediate operator $G_p$ or $G_q$ will change a coefficient of a tropical polynomial by $|p-q|,$ while the points $p$ and $q$ are competing for a single break point (we were depicting such movable points by circles on the previous pictures) and until another break point would reach them.}
\end{figure}

It is easy to verify that the locus of $L(p,q)=1$ has area $\frac 1 4$. It is less trivial to reproduce by hand the locus of $L(p,q)=2$ whose closure consists of two triangles of area $\frac 1 8$, four triangles of area $\frac{1}{16}$ and two triangles of area $\frac{1}{80}$, giving $\frac{21}{40}$ in total.  The loci of $L(p,q)=N\geq2$ are similar to one another, and their areas decrease. Their closures consist of eight triangles (see Figure \ref{lenrelN}) which go in pairs with respect to the symmetry~(3). 

The formula (it is derived simply by summing up the areas of triangles shown on Figure~\ref{lenrelN} -- using equations of their sides, it is easy to find coordinates of all their vertices and compute the areas through determinants) computes the total area as:
\[\frac{3(9N^2-18N+7)}{(3N-1)(3N-2)(3N-4)(3N-5)}\]
which is asymptotically equal to ${\frac 1 3}N^{-2}$ for large $N.$ We conjecture that a similar result holds true for an arbitrary number of points $n.$

\begin{figure}
    \centering
    \includegraphics[width=0.39\textwidth]{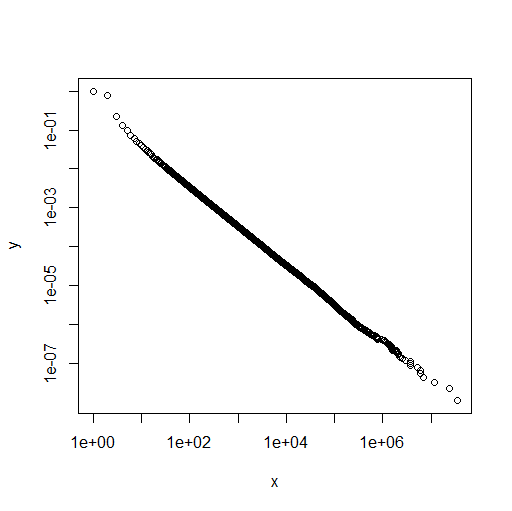}
    \includegraphics[width=0.39\textwidth]{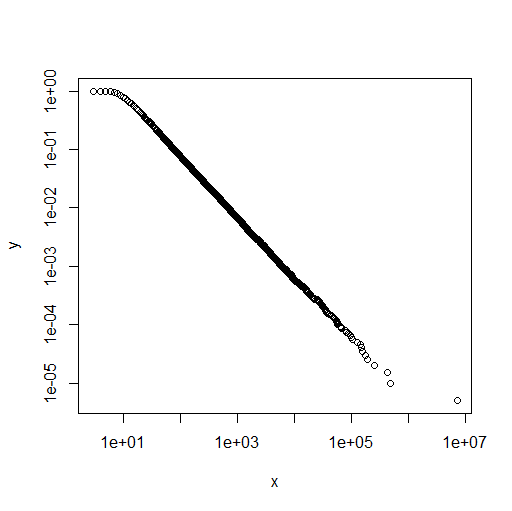}
    \includegraphics[width=0.39\textwidth]{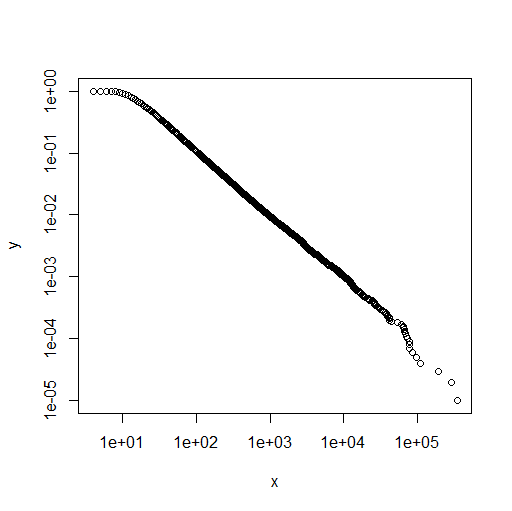}
    \includegraphics[width=0.39\textwidth]{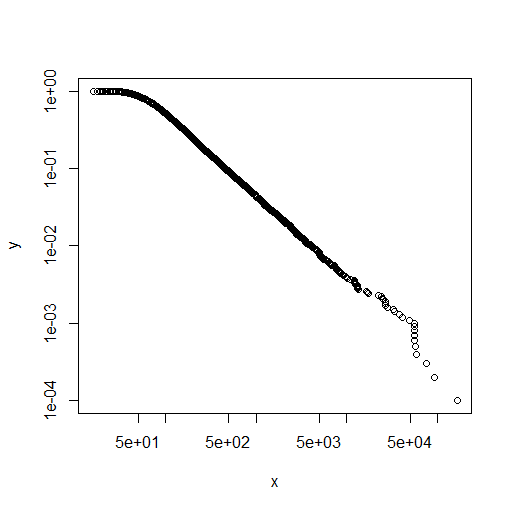}
    \includegraphics[width=0.39\textwidth]{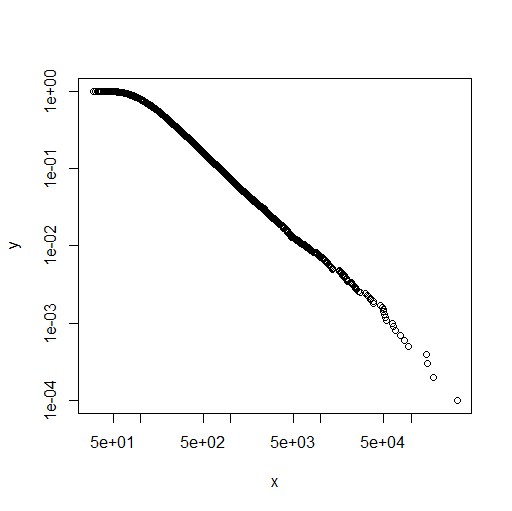}
    \includegraphics[width=0.39\textwidth]{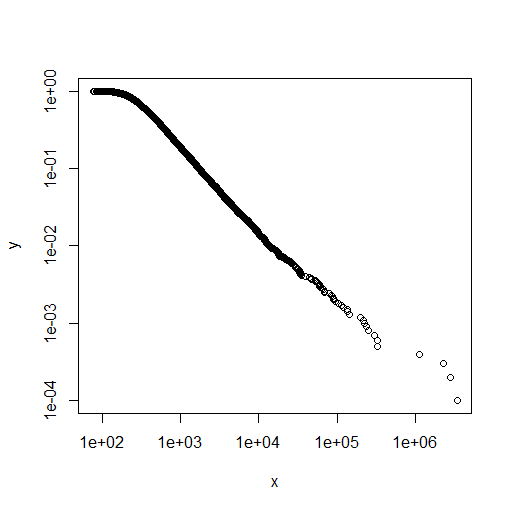}
    \caption{Statistics for the length of relaxation. The y-axis to the $\log$ of the complementary cumulative distribution function and the x-axis corresponds to $\log(L(p_1,\dots,p_n))$. From left to right and from top to bottom: $n=2$ in $10^8$ experiments, $n=7$ in $2\cdot 10^5$ experiments,  $n=8$ in $10^5$ experiments and  $n=16,\ 20$ or $30$ in $10^4$ experiments.}
    \label{plots}
\end{figure}

To justify this, we performed numerical experiments. For a given $n,$ we choose points $p_1,\dots,p_n\in (0,1)$ uniformly at random and gather the statistics of the length of relaxation $L(p_1,\dots,p_n).$ Apart from an anomalous behavior to the left and a noise to the right (due to sporadic appearances of improbably large values), the power-laws are visible (Figure \ref{plots}). 

The computer simulations of relaxations were performed using a program written on OCaml. The data generated through numerous experiments was visualized in R, and the log-log plots in Figure \ref{plots} are obtained using the package poweRlaw \cite{G15}. 

Of course, when looking at the left-hand side of the plots, one can justly object that these are not power-laws in a strict mathematical sense. However, our observable $L$ is conceptually different from those studied in related literature since it can take arbitrarily large values (which is an advantage of the scale-free nature of the model) so we can speak directly about its asymptotic behavior.  We conjecture that for every $n\geq 2$ there exist $\lambda_n<0$ and $c_n>0$ such that \[\text{Measure}(\{{\bf p}\in(0,1)^n:L({\bf p})=N\})\sim c_nN^{\lambda_n}\text{\ as\ }N\rightarrow\infty.\]

Finally, we clarify that spacial observables measuring sizes of avalanches in relaxations are not interesting in the one-dimensional case. For example, we could quantify changes when passing from $G_{p_1,\dots p_n}0_{(0,1)}$ to $G_{p_1,\dots p_n,p_{n+1}}0_{(0,1)}$ by measuring the length of a set $I_{n+1}$ over which these two functions are not equal. This set is easy to find explicitly: let $q_n\in[0,1)$ be the fractional part of $-(p_1+\dots+p_n)$; if $p_{n+1}<q_n,$ then $I_{n+1}=(0,q_n);$ if $p_{n+1}>q_n,$ then $I_{n+1}=(q_n,1)$; and $p_{n+1}\neq q_n$ for generic $p_1,\dots p_{n+1}$. If $p_1,\dots p_{n+1}$ are independent uniform random variables, then $q_n$ is uniform and independent with $p_{n+1}.$ Thus, the distribution of $\text{Length}(I_{n+1})$ doesn't depend on $n\geq 1.$ Its density function is $x\in[0,1]\mapsto 2x.$ 

\subsection*{Acknowledgments}
This work is supported by a grant from The Simons Foundation International (grant no. 992227, IMI-BAS) and by the National Science Fund of Bulgaria, National Scientific Program “VIHREN” (project no. KP-06-DV-7).


\EditInfo{December 14, 2022}{September 10, 2023}{Jacob Mostovoy and Sergei Chmutov}

\end{document}